# Economic Power Capacity Design of Distributed Energy Resources for Reliable Community Microgrids

Chen Yuan[a],*, Guangyi Liu[a], Zhiwei Wang[a], Xi Chen[a], Mahesh S. Illindala[b]

[a]*GEIRI North America, 250 W Tasman Dr., San Jose, CA, 95134, USA*
[b]*The Ohio State University, 2015 Neil Ave., Columbus, OH, 43210, USA*

**Abstract**

Community microgrids are developed within existing power systems by integrating local distributed energy resources (DERs). So power distribution systems can be seamlessly partitioned into community microgrids and end-users could be largely supported when an extreme event happens. However, because of DERs low inertia, their power capacity should be well designed to cover unexpected events and guarantee system reliability. This paper presents a quantitative and qualitative combined methodology for DERs selection, and an economic approach to meet the system reliability requirements. Discrete time Fourier transform (DTFT) and particle swarm optimization (PSO) are employed to obtain the optimal solution, with consideration of load demand and renewable generation uncertainties. In addition, a sensitivity analysis is conducted to show how DERs' capacity design is impacted by counted portion of the forecasted renewable generation.

## 1. Introduction

In recent years, power outages caused by extreme events, like natural disaster and overloading, happened very frequently. To ensure the reliable power supply for electricity consumers, distribution systems could be seamlessly partitioned into community microgrids and each microgrid should be able to operate independently. But one of the existing challenges is how to guarantee the operation reliability of the islanded microgrid. Previous works related to the power capacity design of distributed energy resources (DERs) have been mostly focusing on cost minimization [1], [2]. But few of them ever took generator failure rate and stochastic characteristics of load demand and renewable energy into consideration of system reliability. A study of microgrid generation adequacy indicates that, with a certain generation capacity, more DERs could lead to higher system reliability [3]. The reliability centered generation capacity planning was conducted in [4], but the uncertainty of renewable generation was not included. In this paper, a comprehensive selection of DERs for community microgrids is presented in Section 2. Section 3 studies the impact of planning reserve margin on system reliability. Then, an economic DERs sizing scheme for reliable community microgrids is elaborated in Section 4, with consideration of load uncertainty and renewable's unpredictability. The case study and sensitivity analysis are provided in Section 5, and Section 6 presents a conclusion of this paper.

---

* Corresponding author. Tel.: +1-614-586-3198.
   *E-mail address:* chen.yuan@geirina.net

## 2. Distributed Energy Resources Selection

### 2.1. Levelized cost of energy – Quantitative Assessment

Levelized cost of energy (LCOE) is widely employed by utilities to measure the cost of electricity from a generator. It calculates the generator's annualized cost divided by its estimated yearly energy output. The annualized cost consists of annualized capital cost, operation and maintenance (O&M) cost, fuel cost, and renewable energy tax incentive payback [5], [6]. The O&M cost includes a fixed part and a variable portion. The fixed part is determined by the power rating of each DER, while the variable one is related to the energy output. In this paper, the tax incentive is the production tax credit (PTC) — a federal tax incentive that provides financial support for the development of renewable energy [7].

As a financial tool, LCOE is very valuable for the comparison of various generation units. A low LCOE indicates a low cost of electricity generation. For a conventional power plant, the future fuel price is uncertain and largely depending on external factors, while a renewable energy resource has zero fuel cost, although the initial capital cost is high. Besides, governments have policies to encourage the integration of renewable energy resources, like subsidies, tax incentives, feed-in tariff, net-metering program, renewable portfolio standards, etc. Fig.1 is the comparison of LCOE vs. capacity factor among different DERs.

### 2.2. Qualitative Function Deployment —Qualitative Evaluation

Based on the LCOE assessment, a qualitative evaluation is able to further explore soft indices impact on different types of DERs. This evaluation can be obtained by using qualitative function deployment (QFD). It is a method to convey the customers' voice to engineering evaluation for a product service [8]. In this section, a QFD is employed to examine the relationships — strong (9), medium (3), weak (1) or no relation (0) — between each DER option against environmental factors, customer requirements, and government mandates [9]. Furthermore, it also differentiates between a positive relationship and a negative one. The takeaway from the QFD exercise is able to help better understand which DER(s) need(s) to be considered for community microgrids. As seen in Table 1, biomass generator, natural gas generator, photovoltaic (PV) panel and wind turbine (WT) are more suitable DER options than other choices. Besides, the natural gas generator is very efficient and has ample supplies in the United States.

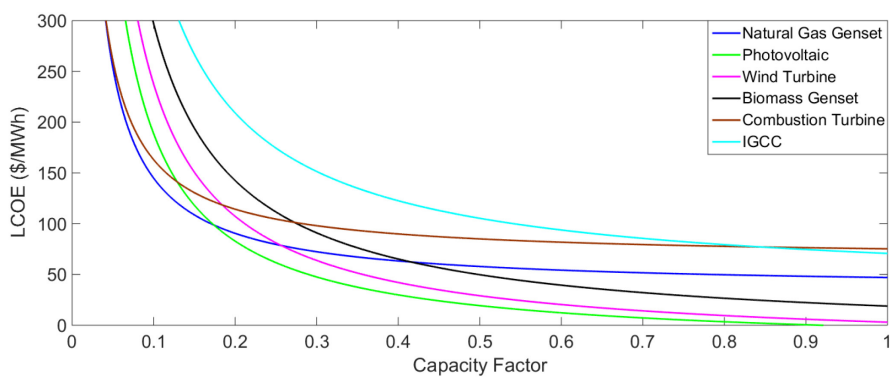

Fig. 1 Curves of LCOE vs. capacity factor for DERs

Table 1. QFD evaluation of DERs for community microgrids development

| DERs Options / Customer Requirements | Importance (1-5) | PV Panel | Wind Turbine | Biomass Generator | Natural Gas Generator | Natural Gas Combustion Turbine | Coal–Fired Power Plant (Base-line) |
|---|---|---|---|---|---|---|---|
| LCOE | 5 | 9 | 9 | 3 | 3 | 1 | 1 |
| $CO_2$ Emission Reduction | 5 | 3 | 3 | 9 | 9 | 1 | 0 |
| Fuel Consumption Savings | 4 | 9 | 9 | 9 | 3 | 1 | 0 |
| Outage Time Reduction | 5 | -3 | -3 | 1 | 3 | 3 | 3 |
| Dispatchability | 4 | -1 | -1 | 1 | 3 | 3 | 1 |
| Equipment Lifetime | 3 | 3 | 3 | 1 | 1 | 1 | 3 |
| Comply with the U.S. DOE Target | 5 | 9 | 9 | 9 | 1 | 1 | 0 |
| **Absolute Target** | | **131** | **131** | **153** | **107** | **49** | **33** |

## 3. Planning Reserve Margin and Its Impact on Reliability

### 3.1. Planning Reserve Margin

Planning reserve margin (PRM) is used to maintain systems reliability. It is a key metric that measures the flexibility to meet load demand and the ability to handle the loss of system components. PRM is usually coupled with probabilistic analysis to identify the resource adequacy and find out whether the generation capacity is large enough to cover peak load demand, loss of one generation unit, and uncertainties from load and renewable resources.

### 3.2. Impact of Planning Reserve Margin on System Reliability

However, system reliability and resource adequacy are not readily observable. For example, it is very difficult to quickly evaluate a system's reliability, like loss of load expectation (LOLE), by simply taking a look at the reserve margin. Based on the probabilistic analysis and Monte Carlo simulation, the general relation between PRM and LOLE is plotted in Fig. 2. Situations like stochastic load changes, renewable intermittency, and one generation loss are studied. Power outages caused by external conditions like extreme weather, grounding fault, and distribution line disconnection are neglected since they can be hardly improved by PRM. Fig. 2 presents the LOLE vs. PRM curves with different proportions of the largest dispatchable generator (PLG). It reflects that, when the largest dispatchable generation unit takes up a larger portion of total generation capacity, the system reliability is worse, needing a larger PRM to achieve the same level of reliability.

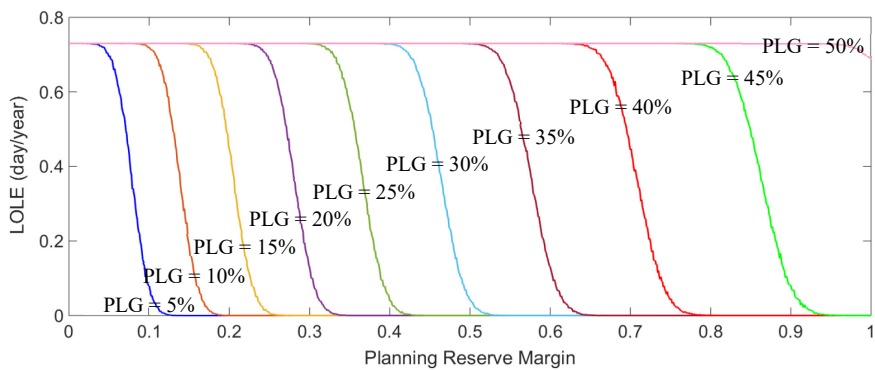

Fig. 2 Curves of LOLE vs. planning reserve margin with different proportions of the largest dispatchable generator

## 4. Power Capacity Design For Dispatchable Generation Units

### 4.1. Problem Formulation

In this paper, renewable energy resources are considered as negative load and sized to meet the customer needs. So the power capacities of renewables like PV and WT, are determined first per customer requirements. The sizing of generators and battery energy storage system (BESS) will be processed then. However, with high penetration of renewable resources, there are extra challenges caused by uncertainties. For example, one question is whether the generators and BESS power capacity design should consider the forecasted generation from renewables or not? If yes, is it better to accommodate a portion of the forecasted renewable generation? These questions will be explored further in Section 4.2 and Section 5.

As described in Section 3, a larger reserve margin will lead to a more reliable system. But it will also result in lower efficiency and higher cost. This is because when the total generation capacity is larger, the operation efficiency could be lower with more reserve margin, resulting in higher capital cost, O&M cost, and fuel cost. So, there is a tradeoff between cost and reliability.

Based on the previous discussion, system reliability is the primary goal of the sizing problem. However, the cost cannot be ignored while carrying out the generators and BESS power capacity design. Therefore, in this paper, the total cost minimization is set as the objective and system reliability requirement is embedded in the constraint. Then, the reliability requirements are satisfied before achieving the cost minimization.

### 4.2. Optimization Algorithm

Based on the problem described in Section 4.1, generators are sized together with BESS to share the net load and provide adequate reserve margin. The net load can be divided into two parts: a) components with large power that varies smoothly over longer duration, and b) small but frequently fluctuating power components. In this way, generators could take the burden of the smooth (i.e., flat) power variation and the BESS can compensate the small and frequent changes.

Discrete time Fourier transform (DTFT) and particle swarm optimization (PSO), which is a population based stochastic optimization technique, are employed to find the optimal power assignment between generators and BESS to minimize annualized cost while satisfying the system reliability requirement. Besides, loss of one generation is covered by PRM to obey "N-1" criterion.

The sizing scheme is explained as follows and illustrated in Fig. 3. With the input of PV and WT power capacities, the net load profile can be achieved based on stochastic models of load and renewable energy resources. The DTFT is applied to obtain the net load frequency spectrum. The spectrum frequency range depends on the sampling rate. Therefore, the time domain net load profile is converted into components in the frequency domain. Then a randomly initialized cut-off frequency divides the net load into two parts. The low-frequency portion of the net load is allocated to generators, while BESS takes care of the high-frequency power components. This helps lower system's capital cost, since the BESS, in terms of power and energy, will avoid being oversized. Once the power share for generators is achieved in the frequency domain, based on the initial cut-off frequency, the power share for generators in the time domain could be obtained by using the method of inverse DTFT. However, the process of inverse DTFT may produce negative values, meaning there is negative power assigned to generators. But the generators cannot absorb power. So the power supply from generator should be adjusted, regulating the negative part to be zero and let BESS take care of it. After making such an adjustment, the BESS power share can be determined by subtracting the generators power allocation from the net load. The next step is to size generators and BESS. The power capacity of generators needs to not only meet the maximum power output but also include a reserve margin to withstand forecast errors and unexpected events.

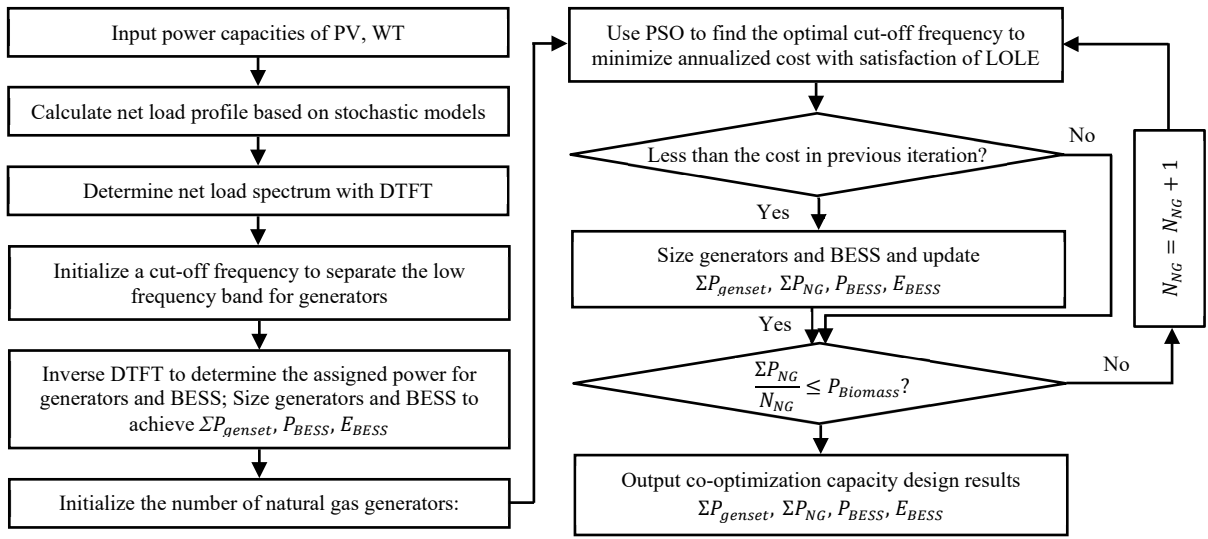

Fig. 3 Flowchart of generators and BESS power capacity design

After the preliminary sizing, the power capacities of generators and BESS are obtained. Power share for natural gas generators equals the subtraction between total power capacity in generators and biomass generator power capacity. However, since the cut-off frequency was initialized randomly, this may not guarantee the optimal power capacity design for both generators and BESS. Therefore, the PSO is employed to find the ideal cut-off frequency and achieve the minimum annualized cost. This is because PSO begins with initialized random solution and searches for optima by updating iterations within the problem space.

In addition, the largest dispatchable generation unit has an impact on system reliability and PRM. Therefore, before the PSO, the power capacity of the largest natural gas generator has to be determined. In Fig. 3, the assumption is that all natural gas generators are identical. So the number of natural gas generators is initialized as one. In each following iteration, the number of natural gas generators is increased by one until the power capacity of the natural gas generator is smaller than the biomass generator power capacity. After all iterations, the optimal solution will be found out.

## 5. Case Study and Sensitivity Analysis

In this section, the proposed strategy for designing generators' and BESS' power capacities is implemented and analyzed for the community microgrid, which has 4 MW peak load, 3 MW PV system, 1 MW wind turbine, and 0.5 MW biomass generator. The load data is generated from a load stochastic model, which is based on two years' historical information provided by the local utility, American Electric Power (AEP) Ohio. The PV and WT output are estimated from their own stochastic models, which are developed using two years' data from sources [10], [11]. In addition, the system requirement of LOLE is set as 1 day in 10 years. The lifetime is assumed for 20 years and the discount rate is 5%.

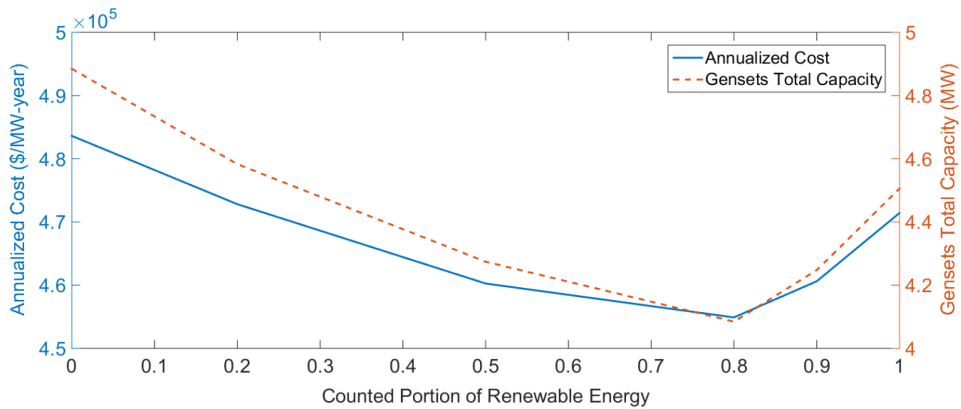

Fig. 4 Tendencies of minimum annualized cost and minimum generators total capacity with counted portion of renewable energy

In this case, all natural gas generators are assumed to be identical. Fig.4 and Table 2 present the comparison of the minimum annualized cost and the minimum total capacity of generators among the six scenarios. It can be easily found that if 80% forecasted renewable generation is counted in the power capacity design, the cost is the minimum. In other words, if we keep 20% renewable energy forecast margin, the net load will be handled in a more economic way. Besides, the situation with consideration of 90% foreseen renewable energy also has a lower minimum annualized cost than the situation of fully considering the forecasted renewable energy.

Table 2. Comparison of Different Scenarios

| Scenarios | Annualized Cost ($/MW-year) | Total Capacity of Generators (MW) | BESS Capacity | | Number of Natural Gas Generators |
|---|---|---|---|---|---|
| | | | Power (MW) | Energy (MWh) | |
| No Renewables | 483,620 | 4.8846 | 0.8507 | 1.1839 | 8 |
| 20% Renewables | 472,810 | 4.5829 | 0.8507 | 1.1581 | 7 |
| 50% Renewables | 460,260 | 4.2742 | 0.8507 | 0.9908 | 7 |
| **80% Renewables** | **454,900** | **4.0846** | **0.8693** | **1.0410** | **6** |
| 90% Renewables | 460,620 | 4.2477 | 0.8693 | 1.0410 | 6 |
| 100% Renewables | 471,420 | 4.5058 | 0.8955 | 1.0308 | 6 |

## 6. Conclusions

This paper presents an economic power capacity design method of generators and BESS for reliable community microgrids. At first, the LCOE based quantitative assessment and QFD based qualitative evaluation are undertaken for various types of DERs to select suitable ones for community microgrids. An economic sizing scheme for generators and BESS is elaborated for ensuring system reliability under uncertainties. The employed optimization methodology is based on DTFT and PSO. In the case study, a sensitivity analysis has been conducted to demonstrate that a small forecast margin of renewable generation could promote cost savings and lower capacity of generators.